\documentclass{amsproc}

\usepackage{amssymb}

\theoremstyle{plain}
\newtheorem{theorem}{Theorem}[section]

\newtheorem{corollary}[theorem]{Corollary}
\newtheorem{definition}[theorem]{Definition}
      
\theoremstyle{remark}

\newcommand{\sgn}{\operatorname{sgn}}
\begin{document}

   \author{W. H. Ko}
   \address{Park Avenue, HKSAR, China}
   \email{wh\_ko@hotmail.com}

   % second author
   
   % title

   \title[Modular Inverse and Reciprocity Formula]{Modular Inverse and Reciprocity Formula}

   % abstract (optional)
   \begin{abstract}
This paper proves a reciprocity formula for modular inverses for non-zero integers and demonstrates some applications of the reciprocity formula in calculating or verifying some modular inverses of specific forms,
including the modular inverse of a Gaussian integer modulo another Gaussian integer.
   \end{abstract}

   % AMS subject classifications (used in AMS journals)

   % AMS keywords (used in AMS journals)

   % acknowledge support, etc
   \thanks{}

   % today's date, or fill in whatever date you prefer
   \date{}

% This ends the top matter information.
% We can now tell LaTeX to display the top matter.

   \maketitle

% Having displayed the top matter, we now proceed to the body of the
% article.

% The body of the article is divided into sections.
% Each section begins with a \section command.

   \section{Introduction}

Modular inverse is one of the basic operations in modular arithmetic,
and it is applied extensively in computer science and telecommunications, particularly, in cryptography.
However, it is also a time-consuming operation implemented in hardware or software compared with other modular arithmetic operations such as addition, subtraction and even multiplication.
Efficient algorithms in calculating modular inverse have been sought after for the past few decades and any improvements will still be welcome.

In implementing an efficient algorithm to calculate the modular inverse of the form $b^{-1} \pmod{2^m}$,
Arazi and Qi \cite{1} have made use of a reciprocity trick which can be translated mathematically to Eq.\eqref{eq:reciprocity}.
As modular inverse is such a fundamental and important operation in mathematics,
it is surprising to find only limited reference to this reciprocity formula in the literature, such as \cite{2},
and none in any elementary textbooks, such as \cite{3,4}.

Although this reciprocity formula can be regarded as a modified form of the linear Diophantine equation Eq.\eqref{eq:lindio},
this formula as written in the form of Eq. \eqref{eq:reciprocity} is able to bring more insight into the properties of modular inverse, as we are going to show in the following discussions.

A recent attempt \cite{6}, which is in Chinese appearing in an IEEE publication, to modernize some ancient Chinese algorithms has introduced a reciprocity formula similar to Eq. \eqref{eq:reciprocity}.  However, its definition of $f_{a,b}$ and $f_{b,a}$ is slightly different from Eq. \eqref{eq:nudef} and thus its reciprocity formula is different.

The classical definition of modular inverse is reproduced in this section, and a new definition of modular inverse is introduced in the next section.
In Section 3, the reciprocity formula and its proof will be provided; and in Section 4, some applications are shown.

\bigskip
\begin{definition}
In modular arithmetic, the classical definition of modular inverse of an integer $a$ modulo $m$ is an integer $x$ such that
\begin{equation} \label{eq:invdef}
ax\equiv 1 \pmod{m}
\end{equation}
The modular inverse is generally denoted as
\begin{equation}
x \equiv a^{-1} \pmod{m}
\end{equation}
\end{definition}

Finding the modular inverse is equivalent to finding the solution of the following linear Diophantine equation, where $a,x,k,m \in\mathbb{Z}$,
\begin{equation}\label{eq:lindio}
a\,x-k\,m=1
\end{equation}

\section{Modular Inverse Definition} 
Slightly different version of modular inverse is used throughout this discussion and it is introduced in this section.
This definition of modular inverse will still satisfy the same congruence equation, Eq.\eqref{eq:invdef}.

Furthermore, in order to have a nicer presentation in the equations containing modular inverse,
new notations for modular operation and modular inverse will be used within this paper.

\begin{definition}
For all $a, m \in\mathbb{Z}$, a modulo m, denoted by $(a)_{m}$, is defined as :
\begin{equation}
(a)_{m}=a-m \bigg\lfloor\dfrac{a}{m}\bigg\rfloor
\end{equation}
Note that \\
\indent \quad \quad
$\begin{cases} 
0 \le (a)_{m} < m & \text{if }m>0\\
m < (a)_{m} \le 0 & \text{if }m<0
\end{cases}$
\end{definition}

\bigskip
In the following discussions, the following definition for modular inverse will be used :

\begin{definition}
Let $a, m \in\mathbb{Z}\setminus \{0\} \text{ and } \gcd(a,m)=1$, modular inverse $a$ modulo $m$, denoted by $(a^{-1})_{m}$, is defined as :
\begin{equation}
(a^{-1})_{m}=x,
\end{equation}
\begin{equation}\label{eq:nudef}
where \begin{cases}
1\le x\le m-1 & \text{if } 1<m \,\&\, a\:x \equiv 1 \pmod{m} \\
m+1\le x \le-1 & \text{if } m<-1 \,\&\, a\,x \equiv 1 \pmod{m} \\
x=\dfrac{1}{2}|m|(\sgn(m)-\sgn(a))+\sgn(a) & \text{if } |m|=1 \\
Undefined & \text{if } a\,m=0 \text{ or } \gcd(a,m) \ne 1
\end{cases}\notag
\end{equation}
\end{definition}

\bigskip

Obviously, for $|m|>1, (a^{-1})_{m}$ satisfies the congruence requirement that
 $$a(a^{-1})_{m} \equiv 1 \pmod{m}$$.

\bigskip

For $|m|=1$, there are two cases.
\begin{description}
\item[Case 1]$m=1\\
\begin{aligned}
(a^{-1})_{m}&=\dfrac{1}{2}|1|(\sgn(1)-\sgn(a))+\sgn(a)=\dfrac{1}{2}(1-\sgn(a))+\sgn(a)\notag\\
&=\begin{cases} 1 & \text{if } a>0 \\ 0 & \text{if } a<0 \end{cases}\notag
\end{aligned}$
\item[Case 2]
$m=-1\\
\begin{aligned}
(a^{-1})_{m}&=\dfrac{1}{2}|-1|(\sgn(-1)-\sgn(a))+\sgn(a)=\dfrac{1}{2}(-1-\sgn(a))+\sgn(a)\notag\\&=\begin{cases} 0 & \text{if } a>0 \\ -1 & \text{if } a<0 \end{cases}\notag.
\end{aligned}$
\end{description}

Hence for $|m|=1, (a^{-1})_{m}$ also satisfies the requirement that $a(a^{-1})_{m} \equiv 1 \pmod{m}$.\\
\indent It is also interesting to note that the modular inverse for $|m|=1$ as defined by Eq.\eqref{eq:nudef} 
is slightly different from the classical definition that $a^{-1} \pmod{m} =0$ for $|m|=1$ and non-zero $a$.

\section{Reciprocity Formula}
The reicprocity relationship between modular inverses seems obvious from the linear Diophantine equation, Eq.\eqref{eq:lindio}.
As a matter of fact, this is the Euclidean algorithm (i.e., iterative division) in disguise.
However, this reciprocity formula is not found in any classic text such as \cite{3}, nor in any modern text, such as \cite{4,5}.
The reciprocity identity first appeared in \cite{2} as Lemma 1 in a format similar to Eq.\eqref{eq:reciprocity}.
However, only positive integers were discussed for specific type of cryptography applications and the conditions that $|m|=1$ was not taken care of.
By the way, on the footnote of p.244 of \cite{2}, it stated that Arazi was "the first to take advantage of this folklore theorem to implement fast modular inversions". 

\begin{theorem}\label{th:reciprocity}\emph{(Reciprocity formula)}
Let $a, b \in\mathbb{Z} \text{ and } \gcd(a,b)=1$, then 

\begin{equation}\label{eq:reciprocity}
a\,(a^{-1})_{b}+b\,(b^{-1})_{a}=1+a \, b. 
\end{equation} 
\end{theorem}

\begin{proof}

\begin{description}
\item[Case 1]$a>1,b>1$\\
Let $U=a(a^{-1})_{b}+b(b^{-1})_{a}$.  Since $U\equiv1\pmod{a}$ and $U\equiv1\pmod{b}$, then $U\equiv1\pmod{a b}.$
That is, $U=1+kab, \text{ where k }\in\mathbb{Z}$. Therefore,\\
$1<a+b \le U=1+k\,a\,b \le a(b-1)+b(a-1)=2a\,b-(a+b)<2a\,b \implies 0<k\,a\,b<2a\,b \implies 0<k<2 \implies k=1$.\\
Therefore $a\,(a^{-1})_{b}+b\,(b^{-1})_{a}=1+a\,b.$

\item[Case 2]$a>1,b<-1$\\
Since $b+1\le (a^{-1})_{b} \le -1 \implies a(b+1) \le a(a^{-1})_{b} \le -a$, and $1 \le (b^{-1})_{a} \le a-1 \implies b(a-1) \le b(b^{-1})_{a} \le b$, therefore, 
$a(b+1) + b(a-1) \le a(a^{-1})_{b} + b(b^{-1})_{a} \le -a+b \implies 2a\,b-(a-b) \le U=1+k\,a\,b \le -(a-b) \implies 2a\,b<2a\,b-(a+1-b) \le k\,a\,b \le -(a+1-b)<0 \implies 0<k<2 \implies k=1$.

\item[Case 3]$a<-1,b>1$\\
Similar to Case 2, and therefore $k=1$.

\item[Case 4]$a<-1,b<-1$\\
Since $b+1 \le (a^{-1})_{b} \le -1 \implies -a \le a(a^{-1})_{b} \le a(b+1)$, and $a+1 \le (b^{-1})_{a} \le -1 \implies -b \le b(b^{-1})_{a} \le b(a+1)$,
therefore $ -a-b \le a(a^{-1})_{b}+b(b^{-1})_{b} \le a(b+1)+b(a+1) \implies -(a+b) \le U=1+k\,a\,b \le 2a\,b+(a+b) \implies 0<-(a+b+1) \le k\,a\,b \le 2a\,b+(a+b+1) < 2a\,b \implies 0<k<2 \implies k=1.$

\item[Case 5]$|a|=1 \text{ or } |b|=1$\\
$a(a^{-1})_{b}+b(b^{-1})_{a}=a(\frac{1}{2}|b|(\sgn(b)-\sgn(a))+\sgn(a))+b(\frac{1}{2}|a|(\sgn(a)-\sgn(b))+\sgn(b))=b(a+\sgn(b))+\sgn(a)(a-a\, b\, \sgn(b))
=(1+a\,b)-(1-a \,\sgn(a))(1-b \,\sgn(b))=(1+a\,b)-(1-|a|)(1-|b|)=1+a\,b.$
\end{description}
\end{proof}

\section{Applications}

\begin{corollary}
If $a,b,k \in\mathbb{Z} \text{ and } \gcd(a,b)=1$, then\\
$((k\,a+b)^{-1})_{a}=\begin{cases} (b^{-1})_{a} & |a|>1\\
(b^{-1})_{a}+\frac{1}{2}(\sgn(k\,a+b)-\sgn(b))) & |a|=1
\end{cases}$
\end{corollary}

\begin{proof}
\begin{description}
\item[Case 1] $|a|>1$\\
$a(a^{-1})_{b}+b(b^{-1})_{a}=1+a\,b \implies (k\,a+b)(b^{-1})_{a}=1+a(b-(a^{-1})_{b}+k(b^{-1})_{a}) \implies
((k\,a+b)^{-1})_{a}=(b^{-1})_{a}$

\item[Case 2] $|a|=1$\\
$((k\,a+b)^{-1})_{a}-(b^{-1})_{a}\\
=\frac{1}{2}|a|(\sgn(a)-\sgn(k\,a+b))+\sgn(k\,a+b)-(\frac{1}{2}|a|(\sgn(a)-\sgn(b))+\sgn(b))\\
=\frac{1}{2}(|a|-2)(\sgn(b)-\sgn(k\,a+b))=\frac{1}{2}(\sgn(k\,a+b)-\sgn(b))$
\end{description}
\end{proof}

\begin{corollary}
If $a,b,k \in\mathbb{Z}, |a| \ne 1$ and $\gcd(a,b)=1$, then
\begin{equation}\label{eq:reduction}
\begin{aligned}
(a^{-1})_{k\,a+b}&=k(a-(b^{-1})_{a})+(a^{-1})_{b}\\
(a^{-1})_{k\,a-b}&=k(b^{-1})_{a}-(b-(a^{-1})_{b})
\end{aligned}
\end{equation}
\end{corollary}

\begin{proof}
\begin{align}
(a^{-1})_{k\,a+b}&=\dfrac{1+a(k\,a+b)-(k\,a+b)((k\,a+b)^{-1})_{a}}{a}\notag\\
&=(k\,a+b)+\dfrac{1-(k\,a+b)(b^{-1})_{a}}{a}\notag\\
&=(k\,a+b)-k(b^{-1})_{a}+\dfrac{1-b(b^{-1})_{a}}{a}\notag\\
&=k(a-(b^{-1})_{a})+(a^{-1})_{b}\notag
\end{align}
This completes the first equation.  And note $((k\,a-b)^{-1})_a=((-b)^{-1})_a=a-(b^{-1})_a$, and the second equation will be obtained.
\end{proof}
It is interesting to note that Eq.\eqref{eq:reduction} fails if the classic definition of $(a^{-1})_b=0 \text{ for }|b|=1$ is used.

On the other hand, although in the proof, it is assumed that $a,b,k \in\mathbb{Z}$, it is a bonus to note that this equation is also valid when $k=i=\sqrt{-1}$.  That is, the modular invesre of $a \pmod{a\,i+b}$ is $(a^{-1})_b+i\,(a-(b^{-1})_a)$.

\bigskip

\begin{corollary}
If $\gcd(a,b)=1$, then
\begin{equation}
\begin{aligned}
((b^{2})^{-1})_{a^{2}}&=(((b(b^{-1})_{a}-2)(b^{-1})_{a})^{2})_{a^{2}}\\
&=((3-2b(b^{-1})_{a})((b^{-1})_{a})^{2})_{a^{2}}
\end{aligned}
\end{equation}
\end{corollary}

\begin{proof}
Let $x=b-(a^{-1})_{b},y=(b(b^{-1})_{a}-2)(b^{-1})_{a}$, and since
$a\:(a^{-1})_{b}+b\:(b^{-1})_{a}=1+a \, b \implies b^{2}y^{2}=1+a^{2}(a^{2}x^{2}-2)x^{2}$, the first part of the equation will become obvious.
Similarly, $b^{2}((3-2b(b^{-1}){a})((b^{-1})_{a})^{2})=1-a^{2}x^{2}(1+2b(b^{-1})_{a})$, and the second part is then also proved.
\end{proof}

\bigskip

For the following corollaries, let $a,b,c,d \in\mathbb{Z}, \gcd(a,b)=\gcd(c,d)=1,\\
u=a\:c+b\:d,v=a\:d-b\:c,s=a^{2}+b^{2},t=c^{2}+d^{2},\\
x_{1}=a(d^{-1})_{c}+b(d-(c^{-1})_{d}),x_{2}=a(c-(d^{-1})_{c})+b(c^{-1})_{d},\\
x_{3}=a(d-(c^{-1})_{d})-b(d^{-1})_{c}),x_{4}=a(c^{-1})_{d}-b(c-(d^{-1})_{c}),\\
y_{1}=c(a-(b^{-1})_{a})+d(a^{-1})_{b},y_{2}=c(b^{-1})_{a}+d(b-(a^{-1})_{b}),\\
y_{3}=c(b-(a^{-1})_{b})-d(b^{-1})_{a},y_{4}=c(a^{-1})_{b}-d(a-(b^{-1})_{a}).
$

\begin{corollary}
If $|u|>1$ and $|v|>1$, then
\begin{equation}\label{eq:co1}
((x_{i})^{-1})_{u}=(y_{i})_{u}, \text{ for }i=1,2.
\end{equation}
\begin{equation}\label{eq:co2}
((x_{i})^{-1})_{v}=(y_{i})_{v}, \text{ for }i=3,4.
\end{equation}
\end{corollary}

\begin{proof}
By direct expansion, and applying Eq.\eqref{eq:reciprocity} where appropriate, it can be shown that\\
$
x_{1}y_{1}=1+u\,((a^{-1})_{b}(d-(c^{-1})_{d})+(d^{-1})_{c}(a-(b^{-1})_{a})),\\
x_{2}y_{2}=1+u\,((b^{-1})_{a}(c-(d^{-1})_{c})+(c^{-1})_{d}(b-(a^{-1})_{b})),\\
x_{3}y_{3}=1-v\,((b^{-1})_{a}(d-(c^{-1})_{d})-(d^{-1})_{c}(b-(a^{-1})_{b})),\\
x_{4}y_{4}=1-v\,((a^{-1})_{b}(c-(d^{-1})_{c})-(c^{-1})_{d}(a-(b^{-1})_{a})).
$

\bigskip
Hence, Eqs.\eqref{eq:co1} and \eqref{eq:co2} are obtained.
\end{proof}

\bigskip
Furthermore, if $a,b,c$ and $d$ are all positive, then\\
$0<c+d\le y_{1}=c(a-(b^{-1})_{a})+d(a^{-1})_{b}\le c(a-1)+d(b-1)=(a\,c+b\,d)-(a+d)<a\,c+b\,d=u$, and hence,\\
\indent $((x_{1})^{-1})_{u}=y_{1}$, if $a,b,c,d>0$

\begin{corollary}
If $|u|>1, |v|>1$ and $\gcd(u,v)=1$, then
\begin{equation}\label{eq:co3}
((a^{2}+b^{2})^{-1})_{u}=((c(a-(b^{-1})_{a})+d(a^{-1})_{b})(v^{-1})_{u})_{u}
\end{equation}
\begin{equation}\label{eq:co4}
((c^{2}+d^{2})^{-1})_{u}=((a(d^{-1})_{c}+b(d-(c^{-1})_{d}))(v^{-1})_{u})_{u}
\end{equation}
\begin{equation}\label{eq:co5}
((a^{2}+b^{2})^{-1})_{v}=((c(a^{-1})_{b}-d(a-(b^{-1})_{a}))(u^{-1})_{v})_{v}
\end{equation}
\begin{equation}\label{eq:co6}
((c^{2}+d^{2})^{-1})_{v}=((a(c^{-1})_{d}-b(c-(d^{-1})_{c}))(u^{-1})_{v})_{v}
\end{equation}
\end{corollary}

\begin{proof}
Let $z_{1}=a(a-(b^{-1})_{a})+b(a^{-1})_{b},z_{2}=c(d^{-1})_{c}+d(d-(c^{-1})_{d}),\\
z_{3}=c(c-(d^{-1})_{c})+d(c^{-1})_{d}.$

\bigskip

Similarly, it can also be shown that\\
$s\,y_{1}=v+u\, z_{1},\\
t\,x_{1}=v+u\, z_{2},\\
s\,y_{4}=u-v\, z_{1},\\
t\,x_{4}=u+v\, z_{3}.$

\bigskip
Hence, Eqs.\eqref{eq:co3}, \eqref{eq:co4}, \eqref{eq:co5} and \eqref{eq:co6} are obtained, respectively.
\end{proof}

\bigskip
As we have already seen in the above example, the reciprocity formula can also be used to verify the modular inverse of Gaussian integers. Specifically the following corollary demonstrates its use in calculating the modular inverses between two Gaussian integers.
\begin{corollary}
Let $a,b,c,d \in\mathbb{Z}\setminus\{0\}$, and $s=a^2+b^2>1, t=c^2+d^2>1$, if $\gcd(s,t)=1$, then
\begin{equation}
(a+i\,b)u+(c+i\,d)v=1+(a+i\,b)(c+i\,d)(a-i\,b)(c-i\,d)\notag
\end{equation}
where $u=(a-i\,b)(s^{-1})_t,v=(c-i\,d)(t^{-1})_s$.
\end{corollary}
That is, $(a+i\,b)^{-1} \pmod{c+i\,d}$ is congruent to $(a-i\,b)(s^{-1})_t \pmod{c+i\,d}$, and $(c+i\,d)^{-1} \pmod{a+i\,b}$ is congruent to $(c-i\,d)(t^{-1})_s \pmod{a+i\,b}$.
\begin{proof}
\begin{flalign}
(a+i\,b)u+(c+i\,d)v&=(a+i\,b)(a-i\,b)(s^{-1})_t+(c+i\,d)(c-i\,d)(s^{-1})_t\notag\\
&=(a^2+b^2)(s^{-1})_t+(c^2+d^2)(s^{-1})_t=s(s^{-1})_t+t(t^{-1})_s\notag\\
&=1+st=1+(a^2+b^2)(c^2+d^2)\notag\\
&=1+(a+i\,b)(a-i\,b)(c+i\,d)(c-i\,d)\notag
\end{flalign}
\end{proof}

\section{$(a^{-1})_1, (a^{-1})_{-1}$}
It may be doubtful why $(a^{-1})_m$ for $|m|=1$ is defined in such a strange way as in Eq.\eqref{eq:nudef},
as the most 'logical' value is 0 since $a \pmod {1}=0$ in the conventional or classical definition of modular arithmetic.
With some examples, we try to arguably convince ourselves that this is actually a better choice than to have $(a^{-1})_1=0$.

First, since $n\equiv0 \pmod{m}$ for all integers $n$ when $|m|=1$, the definition that
\begin{equation}
(a^{-1})_m=\dfrac{1}{2}|m|(\sgn(m)-\sgn(a))+\sgn(a), |m|=1\notag
\end{equation}
does not contradict with the classical definition.

Secondly, as it is shown in Theorem \ref{th:reciprocity}, the reciprocity formula Eq.\eqref{eq:reciprocity} will be valid
for all $a,b\in\mathbb{Z}$, including $|a|=1$ and $|b|=1$ as long as $\gcd(a,b)=1$.

For example, let $m=1,a>0$, with assumption of the classical definition and the reciprocity formula, we have
\begin{center} $a(a^{-1})_1+1\cdot(1^{-1})_a=1+a\cdot1 \implies (1^{-1})_a=1+a$\end{center}
and that is contradictory to the fact that $0<(1^{-1})_a<a$.

On the other hand, with the definition of Eq.\eqref{eq:nudef}, then
\begin{center} $a(a^{-1})_1+1\cdot(1^{-1})_a=1+a\cdot1 \implies a\cdot1+1\cdot1=1+a\cdot1$\end{center}
and that is true.

Take another example with Eq.\eqref{eq:reduction}.  Let $a=7, b=1, k=3$, and assume the classical definition of modular inverse, we have
\begin{flalign}
(7^{-1})_{3\cdot7+1}&=3(7-(1^{-1})_3)+(7^{-1})_1\notag\\
\implies 19=(7^{-1})_{22}&=3(7-1)+0=18\notag\\
\end{flalign}
and that is contradictory, whereas, if $(7^{-1})_1=1$ then both sides are equal.

   \section{Conclusions}

We have proved the reciprocity formula for modular inverse based on a new definition for modular inverse.
The major difference between the new definition and the classic definition arises when $|m|=1$.
The reciprocity formula is also shown to be valid for some Gaussian integers.

The potential use of the reciprocity formula is far more extensive as seen from the few examples in the discussion above, and much further works remain to be carried out.

% This is the end of the last section.

% Finally we create the bibliography or list of references.

% Every LaTeX document must end with \end{document}.


\begin{thebibliography}{99}

\bibitem{1}{O. Arazi and Hairong Qi, "On Calculating multiplicative inverses modulo $2^{m}$",
{\em IEEE Transactions on Computers, } 57(10):1435 - 1438, October 2008. }

\bibitem{2}{Marc Joye and Pascal Paillier, "GCD-Free Algorithms for Computing Modular Inverses",
C.D. Walter et al. (Eds.): {\em CHES 2003, LNCS 2779}, pp.243-253, 2003. Springer-Verlag Berlin Heidelberg 2003.}

\bibitem{3}{G.H. Hardy and E.M.Wright, "An Introduction to the Theory of Numbers", Fifth Edition, Oxford University Press, 1979.}

\bibitem{4}{Song Y. Yan, "Number Theory for Computing", Second Edition, Spring-Verlag Berlin Heidelbert 2000,2002. }

\bibitem{5}{Kenneth Ireland and Michael Rosen, "A Classical Introduction to Modern Number Theory",
Second Edition, 1998, Springer-Verlag.}

\bibitem{6}{Chunrong Zhang and  Jinggang Zhang, "Enrich Strengthen and Expand Shengyubeifenfa Of Chinese Remainder Theorem",
International Conference on Computer Science and Service System (CSSS), 2011, p. 1254 - 1258.}

   \end{thebibliography}
\end{document}